\documentclass[sn-basic]{sn-jnl}% Basic Springer Nature 
\usepackage{array}
\usepackage{graphicx}%
\usepackage{multirow}%
\usepackage{amsmath,amssymb,amsfonts}%
\usepackage{amsthm}%
\usepackage{mathrsfs}%
\usepackage[title]{appendix}%
\usepackage{xcolor}%
\usepackage{textcomp}%
\usepackage{manyfoot}%
\usepackage{subfigure}
\usepackage{natbib} 
\usepackage{booktabs}%
\usepackage{hyperref}
\usepackage{algorithm}%
\usepackage{algorithmicx}%
\usepackage{algpseudocode}%
\usepackage{listings}%
\theoremstyle{thmstyleone}%
\newtheorem{theorem}{Theorem}[section]%  meant for 
\newtheorem{lemma}[theorem]{Lemma}

\theoremstyle{thmstyletwo}%
\newtheorem{remark}{Remark}[section]%

\theoremstyle{thmstylethree}%

\begin{document}

\title[Article Title]{A class of matrix splitting-based fixed-point iteration method for the vertical nonlinear complementarity problem}

\author[1]{\fnm{Yapeng} \sur{Wang}}\email{ypwang0121@163.com}

\author*[1]{\fnm{Xuewen} \sur{Mu}}\email{xwmu@xidian.edu.cn}

\affil[1]{\orgdiv{School of Mathematics and Statistics}, \orgname{Xidian University}, \orgaddress{\city{Xian}, \postcode{710100}, \state{Shanxi}, \country{China}}}

\abstract{In this paper, we propose a class of matrix splitting-based fixed-point iteration (FPI) methods for solving the vertical nonlinear complementarity problem (VNCP). Under appropriate conditions, we present two convergence results obtained using different techniques and estimate the number of iterations required for the FPI method. Additionally, through numerical experiments, we demonstrated that the FPI method surpasses other methods in computational efficiency.}



\maketitle
\section{Introduction} \label{sec1}
In this paper, we consider the vertical nonlinear complementarity problem (VNCP), which seeks to find a $x \in \mathbb{R}^n$ such that
\begin{equation}
    u(x)^Tv(x)=0, u(x) := Ax+\phi(x) \geq 0,v(x):=Bx+\psi(x) \geq 0,
    \label{VNCP}
\end{equation}
where $u(x),v(x) \in \mathbb{R}^n$, $A,B \in \mathbb{R}^{n \times n}$ are given matrices, and $\phi(x),\psi(x) \in \mathbb{R}^n \rightarrow \mathbb{R}^n$ are two given nonlinear functions defined component-wise as $(\phi(x))_i = (\phi(x_i))$ and $(\psi(x))_i = (\psi(x_i))$. This problem was first introduced in \cite{Xie2023}.\\
\indent The VNCP (\ref{VNCP}) has garnered considerable attention in recent years due to its wide-ranging applications in scientific computing and engineering. Notable examples include stability analysis in nonlinear control theory, nonlinear neural networks, and boundary value problems for differential equations with nonlinear terms. For more details, see \cite{Brogliato2003, Hou2022, Goeleven1996}. On the other hand, when $A,B$ are real symmetric matrices and $\phi(x)=p,\psi(x)=q$ are two real vectors, the VNCP \eqref{VNCP} will degenerate into the Vertical Linear Complementarity Problem (VLCP). The study of VLCP has a long history, and there are currently many effective iterative methods available for solving it, such as projected splitting methods \cite{mezzadri2022}, modulus-based matrix splitting methods \cite{he2022}, neural network method \cite{hou2022neural}. \\
\indent Recently, Fu et al. introduced a modulus-based matrix splitting method for solving the VNCP \eqref{VNCP}, see \cite{Fu2024}. This method transformed the VNCP \eqref{VNCP} into an equivalent system of nonlinear equations
\begin{equation}
    (A+\Omega B) x=|(A-\Omega B) x+(\phi(x)-\Omega \psi(x))|-\phi(x)-\Omega \psi(x),
\end{equation}
where $\Omega$ is a positive diagonal matrix. Subsequently, by splitting the matrices as $A = M_A-N_A$ and $B=M_B-N_B$, an iterative framework is established
\begin{equation}
\left(M_A+\Omega M_B\right) x^{k+1}=\left(N_A+\Omega N_B\right) x^k+|(A-\Omega B) x^k+(\phi(x^k)-\Omega \psi(x^k))|-\phi(x^k)-\Omega \psi(x^k),
\label{VNCP_Fu}
\end{equation}
where $M_A+\Omega M_B$ is nonsingular. Additionally, Fu et al. presented the necessary convergence conditions and validated the effectiveness of the modulus-based matrix splitting method \eqref{VNCP_Fu} through some numerical examples.\\
\indent Due to its simplicity, the fixed-point iteration (FPI) method is widely used in numerical computation and engineering applications, including solving linear systems \cite{Menard2008,Do2008,Jakoveti2021}, absolute value equations \cite{Yu2022,Ali2023,Li2023,Liu2024,Ali2022}, and both linear and nonlinear complementarity problems \cite{Noor1988,Dai2020,Shi2016,Fang2021,Wang2023}. For further details, refer to the cited literature. To the best of our knowledge, no research has yet applied the fixed point iteration method to solve the VNCP \eqref{VNCP}, motivating our study. However, due to the nonlinear components, the FPI method cannot be directly applied to solve the VNCP \eqref{VNCP}. Therefore, in this paper, we aim to establish a fixed-point iterative framework for solving the VNCP. Our contributions are outlined as follows:\\
\indent 1. By introducing an auxiliary variable
$$y=|(A-\Omega B) x+\left(\phi\left(x\right)-\Omega \psi\left(x\right)\right)|,$$ and splitting the coefficient matrix $A+\Omega B$, we construct a class of matrix splitting-based fixed-point iteration (FPI) framework for solving the VNCP\eqref{VNCP}. The FPI method does not require the computation of the inverse of matrix $A+\Omega B$, making it computationally efficient with fast iteration.\\
\indent 2. We provide two convergence results for the FPI method under appropriate assumptions, thereby determining the range of value for the parameter. Additionally, the number of iterations required for the FPI method is also estimated.\\
\indent 3. We present several numerical examples to demonstrate the effectiveness and practical feasibility of our method.\\
\indent This paper is organized into five main sections. In Section \ref{sec2}, we introduce a fixed-point iteration framework to solve the VNCP \eqref{VNCP}. In Section \ref{sec3}, we establish the convergence conditions for the FPI method. The results of numerical experiments and the conclusions of this paper are presented in Sections \ref{sec4} and \ref{sec5}, respectively.

\section{The fixed point iteration method for the VNCP \eqref{VNCP}} \label{sec2}
\indent First, we present some concepts and lemmas that will be used in the subsequent discussion. \\
\indent Let $\mathbb{R}^{n \times n}$ be the set of $n \times n$ real matrices, with $\mathbb{R}^n = \mathbb{R}^{n \times 1}$. The symbol $I_n$ denotes the $n \times n$ identity matrix. Define $|x|$ to represent taking the absolute value of each component of vector $x\in \mathbb{R}^n $: $|x| = (|x_1,\cdots,|x_n|)^{T}$. For any $A \in \mathbb{R}^{n \times n}$, let $\| A \|$ and $\rho(A)$ denote the spectral norm and spectral radius of $A$, defined respectively as $\|A\|:=\max \left\{\|A x\|_2: x \in \mathbb{R}^n,\|x\|_2=1\right\}$ and $\rho(A) := \max\limits_{1 \leq i \leq n} \left\{|\lambda_i|\right\} $, where $\lambda_i$ is an eigenvalue of $A$.
\begin{lemma} \emph{\cite{Fu2024}}
Suppose that $\Omega$ is a positive diagonal matrix, the VNCP \eqref{VNCP} is equivalent to the following nonlinear equations: 
\begin{equation}
(A+\Omega B) x=|(A-\Omega B) x+(\phi(x)-\Omega \psi(x))|-\phi(x)-\Omega \psi(x) .
\label{VNCP_equ}
\end{equation}  
\end{lemma} \label{lemma2.1}
\begin{lemma} \emph{\cite{Bai2021}}
    For any vectors $x_1 \in \mathbb{R}^n$ and $x_2 \in \mathbb{R}^n$, the following conclusions hold:\\
    (1) $\| |x_1|-|x_2| \| \leq \|x_1-x_2\|;$
    (2) if $x_1 \leq x_2$, then $\| x_1 \| \leq \|x_2\|;$
    (3) if $x_1 \leq x_2$ and $P$ is a nonnegative matrix, then $Px_1 \leq Px_2$.
    \label{lemma2.2}
\end{lemma}
\begin{lemma} \emph{\cite{Young2014}}
    For real numbers $m$ and $n$, the roots of the quadratic equation $x^2-mx+n=0$ are less than 1 in moduls if and only if $\lvert n \rvert < 1$ and $\lvert m \rvert < 1+n.$
    \label{lemma2.3}
\end{lemma}

Next, we consider designing a fixed point iteration method to solve the nonlinear equations \eqref{VNCP_equ}. \\
\indent Let $A+\Omega B = M-N$, where $M$ is invertible, and define $y=|(A-\Omega B) x+(\phi(x)-\Omega \psi(x))|$, then we have 
\begin{equation}
\left\{\begin{array}{l}
x^*=M^{-1} [Nx^*+y^*-\phi(x^*)-\Omega \psi(x^*)], \\
y^*=(1-\tau)y^*+\tau |(A-\Omega B)x^*+\phi(x^*)-\Omega \psi(x^*)|,
\end{array}\right. \label{VNCP_exact}
\end{equation}
where $\tau$ is a positive constant and $x^*$ denotes the exact solution of the VNCP \eqref{VNCP}. \\
\indent From the above, it can be seen that a class of matrix splitting-based fixed point iteration (FPI) framework for VNCP can be established
\begin{equation}
\left\{\begin{array}{l}
x^{k+1}=M^{-1} [Nx^k+y^k-\phi(x^k)-\Omega \psi(x^k)], \\
y^{k+1}=(1-\tau)y^k+\tau |(A-\Omega B)x^{k+1}+\phi(x^{k+1})-\Omega \psi(x^{k+1})|.
\end{array}\right. \label{VNCP_iterative}
\end{equation}
\indent Clearly, $A+\Omega B = M-N$ provides various forms for the FPI method \eqref{VNCP_iterative}, including different iteration methods in special cases. We decompose $A+\Omega B = D-L-U$, where $D$, $-L$, and $-U$ represent the diagonal matrix, the strictly lower triangular matrix, and the strictly upper triangular matrix of $A+\Omega B$, respectively. Then
\begin{itemize}
    \item Let $M=D$, $N=L+U$, the FPI method \eqref{VNCP_iterative} gives the fixed point-based Jacobi (FPI-J) iteration method.
    \item Let $M=D-L$, $N=U$, the FPI method \eqref{VNCP_iterative} gives the fixed point-based Gauss-Seidel (FPI-GS) iteration method.
    \item Let $M=\frac{1}{\alpha}D-L$, $N=(\frac{1}{\alpha}-1)D+U$, the FPI method \eqref{VNCP_iterative} gives the fixed point-based successive overrelaxation (FPI-SOR) iteration method, where $ 0<\alpha<2.$
\end{itemize}

\indent We note that the FPI method \eqref{VNCP_iterative} fully considers the splitting of the coefficient matrix $A+\Omega B$. As a result, it provides multiple solutions for solving more complex problems in practical computations.
\section{Convergence analysis}\label{sec3}
In this section, we analyze the convergence of the FPI method \eqref{VNCP_iterative} to determine the range of value for the parameter $\tau$. To this end, we give some notations relevant to the FPI method \eqref{VNCP_iterative}:
\begin{equation*}
    \varepsilon_x^k = x^k-x^*, \varepsilon_y^k = y^k-y^*,  E^{(k+1)}  =\left[\begin{array}{l}
\| \varepsilon_x^{k+1} \| \\
\| \varepsilon_y^{k+1} \|
\end{array}\right].
\end{equation*}
Assume there exist positive constants $L_1$ and $L_2$ such that 
    \begin{equation*}
    \| \phi(u)-\phi(v) \| \leq L_1 \| u-v\|, \, \| \psi(u)-\psi(v) \| \leq L_2 \| u-v\|, \, \text{for any vectors } u,v \in \mathbb{R}^n.
    \end{equation*}
Next, we present the convergence result of the FPI method \eqref{VNCP_iterative}. 

\begin{theorem}\label{theo3.1}
    Let $A+\Omega B = M-N$ be the matrix splitting form of $A+\Omega B$, where $M$ is non-singular.
    Denote $\alpha = \| M^{-1} \|, \beta = \| N\| + L_1 +L_2 \| \Omega \|, \gamma = \| A-\Omega B\|+L_1+L_2 \| \Omega \|.$ \\
    If 
    \begin{equation}
    \alpha(\beta+\gamma)<1 \quad and 
    \quad 0 < \tau < \frac{2(1-\alpha \beta)}{1-\alpha \beta +\alpha \gamma},\label{conver_condi}
    \end{equation}
    then the sequence $\{ x^k\}_{k=0}^{+\infty}$ obtained by the FPI method \eqref{VNCP_iterative} converges to the unique solution $x^*$ of the VNCP \eqref{VNCP}.
\end{theorem}
\textbf{Proof.} It follows from \eqref{VNCP_exact} and \eqref{VNCP_iterative}, we have
\begin{equation}
    \begin{aligned}
    \varepsilon_{x}^{k+1} = M^{-1}[N\varepsilon_x^k+\varepsilon_y^k+(\phi(x^*)-\phi(x^k))+\Omega(\psi(x^*)-\psi(x^k))],
    \label{wucha_x}
    \end{aligned}
\end{equation}
and 
\begin{equation}
    \begin{aligned}
    \varepsilon_{y}^{k+1} = (1-\tau)\varepsilon_y^k+\tau[|(A-\Omega B)x^{k+1} + \phi(x^{k+1})-\Omega \psi(x^{k+1})| - |(A-\Omega B)x^{*} + \phi(x^*)-\Omega \psi(x^{*})|].\\
    \label{wucha_y}
    \end{aligned}
\end{equation}
From \eqref{wucha_x}, we obtain
\begin{equation}
\begin{aligned}
    \| \varepsilon_x^{k+1}\| &\leq \|M^{-1} \|\cdot (\left\|N\|+L_1+L_2\|\Omega \|\right)\cdot \| \varepsilon_x^{k}\|+\|M^{-1}\| \cdot \|\varepsilon_y^k\|\\
    &=\alpha \beta \cdot \| \varepsilon_x^k\|+\alpha \cdot \|\varepsilon_y^k\|.
    \label{wucha_norm_x}
    \end{aligned}
\end{equation}
From \eqref{wucha_y} and Lemma \ref{lemma2.2}, we obtain
\begin{equation}
    \begin{aligned}
    \| \varepsilon_y^{k+1}\| &\leq |1-\tau| \cdot \|\varepsilon_y^k\| +\tau \| (A-\Omega B)(x^{k+1}-x^*) + \phi(x^{k+1})-\phi(x^*)+\Omega (\psi(x^*)-\psi(x^{k+1})) \|\\
    &\leq |1-\tau| \cdot \|\varepsilon_y^k\| + \tau (\|A-\Omega B\|+L_1+L_2\| \Omega \|)\cdot \|\varepsilon_x^{k+1}\| \\
    & = |1-\tau| \cdot \|\varepsilon_y^k\| + \tau \gamma \cdot \|\varepsilon_x^{k+1}\|.
    \label{wucha_norm_y}
    \end{aligned}
\end{equation}
Hence, from \eqref{wucha_norm_x} and \eqref{wucha_norm_y}, we have
\begin{equation}
\left[\begin{array}{cc}
1 & 0 \\
-\tau \gamma & 1
\end{array}\right] 
\begin{bmatrix} \left\|\varepsilon_{x}^{x+1}\right\| \\ \left\|\varepsilon_{y}^{k+1}\right\| \end{bmatrix} 
\leq 
\left[\begin{array}{cc}
\alpha \beta & \alpha \\
0 & |1-\tau|
\end{array}\right] 
\begin{bmatrix} \left\|\varepsilon_x^k\right\| \\ \left\|\varepsilon_y^k\right\| \end{bmatrix} .
\label{wucha_xy}
\end{equation}
\indent Let \begin{equation*}
P = \left[\begin{array}{ll}
1 & 0  \\
\tau \gamma & 1  \\
\end{array}\right],
\end{equation*} multiplying \eqref{wucha_xy} on the left by the nonnegative matrix $P$. We get
\begin{equation}
\left[\begin{array}{l}
\| \varepsilon_x^{k+1} \| \\
\| \varepsilon_y^{k+1} \|
\end{array}\right]
\leq 
\left[\begin{array}{cc}
\alpha \beta & \alpha \\
\tau \gamma \alpha \beta & \tau \gamma \alpha+ |1-\tau|
\end{array}\right] 
\left[\begin{array}{l}
\| \varepsilon_x^{k} \| \\
\| \varepsilon_y^{k} \|
\end{array}\right] .
\label{wucha_xy}
\end{equation}
Then, we have
\begin{equation*}
    E^{(k+1)}  \leq T  E^{(k)}.
\end{equation*}
\indent Now, we consider the value of parameter $\tau$ to ensure that $\rho(T) < 1$. Assume that $\lambda$ is an eigenvalue of $T$, we get
\begin{equation*}
    \lambda^2 - (\alpha \beta +\tau \gamma \alpha +|1-\tau|)\lambda +\alpha \beta |1-\tau| = 0.
\end{equation*}
According to Lemma \ref{lemma2.3}, $|\lambda| < 1$ if and only if
\begin{equation}
\left\{\begin{array}{l}
\alpha \beta + \tau \gamma \alpha + |1-\tau| <1+\alpha \beta |1-\tau|, \\
\alpha \beta |1-\tau| < 1.
\end{array}\right.  \label{lambda_gen}
\end{equation}
Next, we will analyze various cases.
\begin{itemize}
    \item Case 1: When $0<\tau \leq 1$, \eqref{lambda_gen} can be reduced to 
    \begin{equation}
    \left\{\begin{array}{l}
    \tau (\alpha \gamma + \alpha \beta ) < \tau, \\
    1 - \frac{1}{\alpha \beta} <\tau.
    \end{array}\right.  
    \label{gen}
    \end{equation}
    Since $0<\tau \leq 1$ and $\alpha \gamma + \alpha \beta <1$, then \eqref{gen} holds naturally.

    \item Case 2: When $1<\tau$, \eqref{lambda_gen} can be reduced to 
    \begin{equation*}
    \left\{\begin{array}{l}
    \tau < \frac{2(1-\alpha \beta)}{1-\alpha \beta +\alpha \gamma}, \\
    \tau < \frac{1}{\alpha \beta}+1.
    \end{array}\right.  \label{gen_01}
    \end{equation*}
    Since $\alpha \beta + \alpha \gamma <1$, we obtain 
    \begin{equation*}
        1< \frac{2(1-\alpha \beta)}{1-\alpha \beta +\alpha \gamma}  < \frac{1}{\alpha \beta} + 1.
    \end{equation*}Hence, in this case, we get
    \begin{equation*}
    1< \tau < \frac{2(1-\alpha \beta)}{1-\alpha \beta +\alpha \gamma}
    \end{equation*}
\end{itemize}
In summary, we obtain $\rho(T) < 1$ if and only if
\begin{equation*}
    \left\{\begin{array}{l}
    \alpha(\beta+\gamma)<1, \\
    0 < \tau < \frac{2(1-\alpha \beta)}{1-\alpha \beta +\alpha \gamma}.
    \end{array}\right.  \label{gen_01}
\end{equation*}
\indent This implies that $ E^{k+1}  = 0 (k \rightarrow \infty )$ and the sequence $\{ x^k\}_{k=0}^{+\infty}$ generated by the FPI method \eqref{VNCP_iterative}  converges to the unique solution of the VNCP \eqref{VNCP}. The proof is completed.\\
\indent Next, by using an error estimate with a new weighted norm, we derive an alternative scheme that achieves convergence result of the FPI method \eqref{VNCP_iterative}.
\begin{theorem}\label{theo3.2}
    Let the definitions of $\alpha$, $\beta$ and $\gamma$ be the same as in Theorem \ref{theo3.1}. Denote
    \begin{equation*}
      E_{\gamma}^{(k+1)}  =\left[\begin{array}{l}
    \gamma \| \varepsilon_x^{k+1} \| \\
           \| \varepsilon_y^{k+1} \|
    \end{array}\right].
    \end{equation*}
    Then we have 
    \begin{equation}
    \| E_{\gamma}^{(k+1)} \|_{\infty} \leq \| T_{\gamma} \|_{\infty}  \| E_{\gamma}^{(k)} \|_{\infty},
    \label{new_error}
    \end{equation}
    where $$T_{\gamma}=
    \left[\begin{array}{cc}
\alpha \beta & \alpha \gamma\\
\tau \alpha \beta & \tau \alpha \gamma + |1-\tau|
\end{array}\right] .$$
Furthermore, $\| T_{\gamma} \|_{\infty}   < 1$ if and only if 
\begin{equation}
    \left\{\begin{array}{l}
    \alpha(\beta+\gamma)<1, \\
    0 < \tau < \frac{2}{\alpha(\beta+\gamma)+1}.
    \end{array}\right.  \label{new_convergence}
\end{equation}
Then if \eqref{new_convergence} holds, the sequence $\{ x^k\}_{k=0}^{+\infty}$ obtained by the FPI method \eqref{VNCP_iterative} converges to the unique solution $x^{*}$ of the VNCP \eqref{VNCP}.
\end{theorem}

\textbf{Proof.} Let
    \begin{equation*}
        Q = \left[\begin{array}{ll}
        \gamma & 0  \\
        0  & 1  \\
        \end{array}\right]>0,
    \end{equation*}
according to Lemma \ref{lemma2.2}, we multiply \eqref{wucha_xy} on the left by the nonnegative matrix $Q$ to get
\begin{equation} \left[\begin{array}{l}
    \gamma \| \varepsilon_x^{k+1} \| \\
           \| \varepsilon_y^{k+1} \|
    \end{array}\right]
\leq 
\left[\begin{array}{cc}
\alpha \beta & \alpha \gamma  \\
\tau \alpha \beta & \tau \alpha \gamma + |1-\tau|
\end{array}\right] 
\begin{bmatrix} \gamma\left \|\varepsilon_x^k\right\| \\ \left\|\varepsilon_y^k\right\| \end{bmatrix},
\end{equation}
this indicates that
    \begin{equation}
     E_{\gamma}^{(k+1)} \leq  T_{\gamma}   E_{\gamma}^{(k)}.
    \end{equation}
Hence, it can be concluded that \eqref{new_error} holds.\\
\indent Since $\| T_{\gamma} \|_{\infty} = \max\{\alpha(\beta+\gamma),\tau \alpha (\beta+\gamma) + |1-\tau| \}$. In order for $\| T_{\gamma} \|_{\infty}$ to be less than $1$, we obtain the following inequalities:
\begin{equation}
    \left\{\begin{array}{l}
    \alpha(\beta+\gamma)<1, \\
    \tau \alpha (\beta+\gamma) + |1-\tau| < 1.
    \end{array}\right.  
\end{equation}
Therefore, we have
\begin{equation*}
\left\{\begin{array} { l } 
{ \alpha(\beta+\gamma) < 1 } \\
{ 1-\tau \alpha (\beta+\gamma)>0 } \\
{ \tau \alpha(\beta+\gamma) -1 <1-\tau<1-\tau \alpha (\beta+\gamma) }
\end{array} \Leftrightarrow \left\{\begin{array}{l}
{ \alpha(\beta+\gamma) < 1 } \\
{ 0<\tau< \frac{1}{\alpha(\beta+\gamma)} } \\
{ 0<\tau< \frac{2}{\alpha(\beta+\gamma)+1} }
\end{array}\right.\right.
\end{equation*}
\begin{equation*}
\Leftrightarrow\left\{\begin{array} { l } 
{ \alpha(\beta+\gamma) < 1 } \\
{ 0<\tau< \frac{2}{\alpha(\beta+\gamma)+1}}.
\end{array} \right.
\end{equation*}
From \eqref{new_error}, we have
\begin{equation*}
0 \leq\|E^{(k)}_{\gamma}\|_{\infty} \leq\|T_{\gamma} \|_{\infty} \cdot\|E^{(k-1)}_{\gamma}\|_{\infty} \leq \cdots \leq\|T_{\gamma}\|_{\infty}^k \cdot\|E^{(0)}_{\gamma} \|_{\infty}.
\end{equation*}
\indent Hence, if \eqref{new_convergence} holds, we obtain $\lim\limits_{k \rightarrow \infty}\|E^{(k)}_{\gamma}\|_{\infty}=0$. Since $\|E^{(k)}_{\gamma}\|_{\infty} = \max\{ \gamma \| \varepsilon_x^{k} \| ,\| \varepsilon_y^{k} \|\}$, it follows that 
$$
\lim _{k \rightarrow \infty}\left\|\varepsilon_x^k\right\|=0 \quad \text { and } \quad \lim _{k \rightarrow \infty}\left\|\varepsilon_y^k\right\|=0.
$$
From above that, we conclude that the sequence $\{ x^k\}_{k=0}^{+\infty}$ obtained by the FPI method \eqref{VNCP_iterative} converges to the unique solution $x^*$ of the VNCP \eqref{VNCP} under the condition \eqref{new_convergence}.

\begin{remark}
\end{remark}
Comparing the results obtained from Theorem \ref{theo3.1} and Theorem \ref{theo3.2}, since $\alpha (\beta + \gamma) <1$, we have 
$$ \frac{2}{\alpha(\beta+\gamma)+1} < \frac{2(1-\alpha \beta)}{1-\alpha \beta +\alpha \gamma}.$$ This indicates that the result of Theorem \ref{theo3.2} is a stricter form of Theorem \ref{theo3.1}, which may yield better performance in practical computations.\\
\indent Next, we attempt to use the result from Theorem \ref{theo3.2} to estimate the iteration step of the FPI method.
\begin{theorem}
    If \eqref{new_convergence} holds, then the estimate of iteration step $k$ satisfies
    \begin{equation*}
    k>\log_{\xi} c,
    \end{equation*}
    where $c$ and $\xi$ are positive constants.
\end{theorem}
\textbf{Proof.} Based on the results from Theorem \ref{theo3.2}, we have 
\begin{equation*}
    \left\{\begin{array}{l}
    \| T_{\gamma} \| < 1, \\
    0 \leq \| E_{\gamma}^{(k)} \|_{\infty} \leq \| T_{\gamma} \|_{\infty} \cdot \| E_{\gamma}^{(k-1)} \|_{\infty} \leq \cdots \leq \| T_{\gamma} \|_{\infty}^k \cdot \| E_{\gamma}^{(0)} \|_{\infty}.
    \end{array}\right.  
\end{equation*}
Then we assume that
\begin{equation*}
\|E^{(k)}_{\gamma}\|_{\infty} \leq\|T_{\gamma}\|_{\infty}^k \cdot\|E^{(0)}_{\gamma} \|_{\infty} < \delta,
\end{equation*}
thus, we have
\begin{equation*}
\|T_{\gamma}\|_{\infty}^k < \frac{\delta}{\|E^{(0)}_{\gamma} \|_{\infty}}.
\end{equation*} Here, we set $c=\frac{\delta}{\|E^{(0)}_{\gamma} \|_{\infty}}, 0<\|T_{\gamma}\|_{\infty}<\xi$, we obtain
\begin{equation*}
    k>\frac{\log c}{\log (\|T_{\gamma}\|_{\infty})}>\frac{\log c}{\log \xi} = \log_{\xi} c.
\end{equation*}
\indent The proof is completed.

\section{Numerical experiments}\label{sec4}
\begin{figure}[!h]
    \centering
	\subfigure[FPI-J]{
			\includegraphics[scale=0.27]{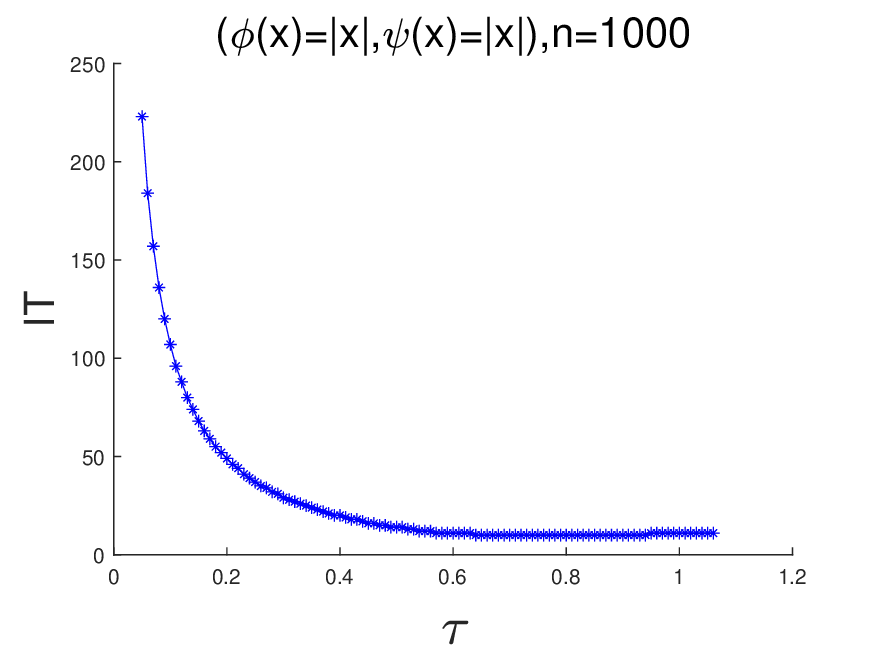} \label{1}
		}
	\subfigure[FPI-GS]{
			\includegraphics[scale=0.27]{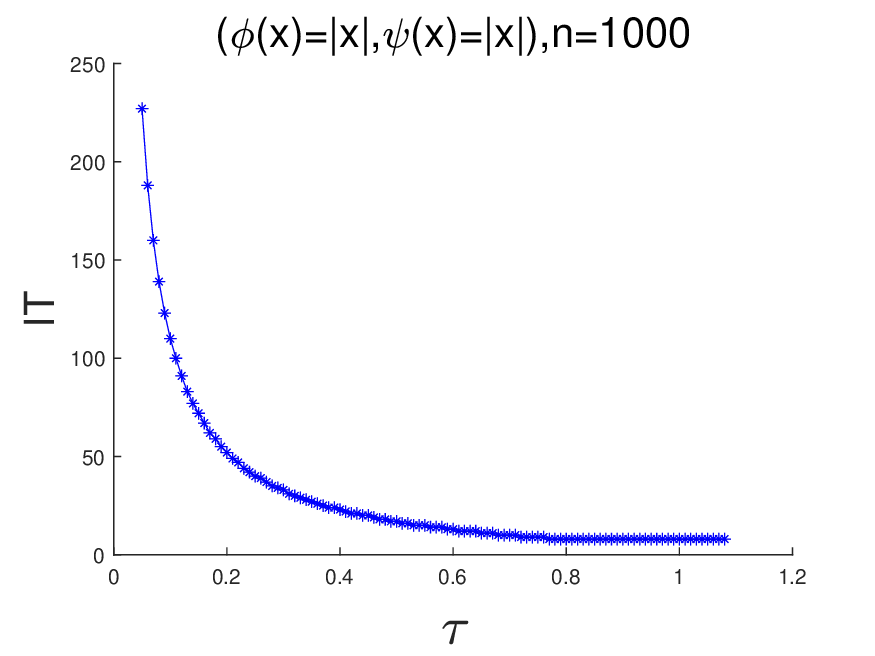} \label{2} 
		}
        \subfigure[FPI-SOR]{
			\includegraphics[scale=0.27]{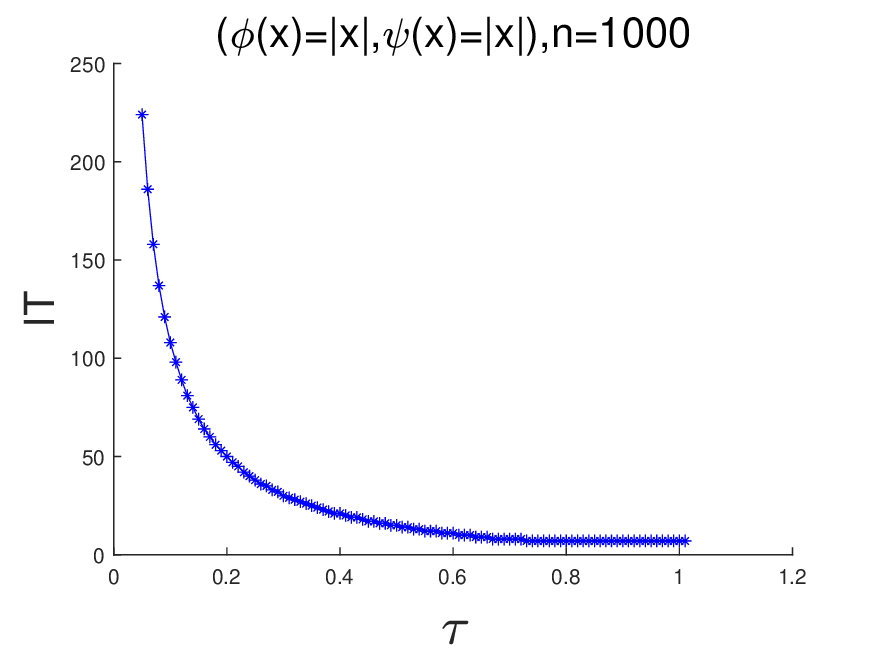} \label{3} 
		}
	\caption{The IT of Example 4.1 for FPI method \eqref{VNCP_iterative}.} \label{fig1}
\end{figure}

In this section, we use several numerical examples to validate the effectiveness and superiority of the FPI method \eqref{VNCP_iterative}. For the different splitting forms of $A+\Omega B$, we define them as FPI-J, FPI-GS, and FPI-SOR, respectively. Two criteria are used to assess the performance of a method: ‘‘IT’’ denotes the number of iteration steps, and ‘‘Time(s)’’ denotes the elapsed CPU time. In our experiments, we use \textbf{RES} $< 10^{-8}$ as the stopping criteria, which is defined as
\begin{equation*}
    \textbf{RES} = |A x^k+\phi(x^k)|^T|B x^k+\psi(x^k)|,
\end{equation*}
the initial vector for all methods is defined as  $x^{(0)} = (1,1,...,1,1) \in \mathbb{R}^n$, and $\Omega = 5I$. All experiments are implemented on a personal computer using MATLAB 2023a (Intel(R) Core i7-12650H, 16GB RAM).\\
\indent We adopt the following strategy to select the iteration parameter $\tau$: first, calculate the range of $\tau$ based on the convergence results derived from Theorem \ref{theo3.1}. Then, plot the relationship between the number of iterations (IT) and $\tau$, allowing us to identify the optimal parameter visually. To provide a more intuitive demonstration, in Figure \ref{fig1}, we describe the iterations (IT) of our methods for solving Example 4.1 with $\phi(x)=\psi(x)=|x|$, respectively.

\textbf{Example 4.1}
Let 
\begin{equation*}
A=\left[\begin{array}{ccccc}
8 & -1 & & & \\
-1 & 8 & -1 & & \\
& \ddots & \ddots & \ddots & \\
& & -1 & 8 & -1 \\
& & & -1 & 8
\end{array}\right] \in \mathbb{R}^{n \times n}, B=\left[\begin{array}{ccccc}
4 & -1 & & & \\
-1 & 4 & -1 & & \\
& \ddots & \ddots & \ddots & \\
& & -1 & 4 & -1 \\
& & & -1 & 4
\end{array}\right]\in \mathbb{R}^{n \times n} \text {. }
\end{equation*}

\textbf{Example 4.2}
Let $A = \hat{A}+\mu_1 I$, $B = \hat{B}+\mu_2 I$, where $\mu_1,\mu_2 \in \mathbb{R}$,
\begin{equation*}
\hat{A}=\left[\begin{array}{ccccc}
S & -I & & & \\
-I & S & -I & & \\
& \ddots & \ddots & \ddots & \\
& & -I & S & -I \\
& & & -I & S
\end{array}\right]\in \mathbb{R}^{n \times n}, \hat{B}=\left[\begin{array}{ccccc}
S &  & & & \\
 & S &  & & \\
&  & \ddots &  & \\
& &  & S &  \\
& & &  & S
\end{array}\right]\in \mathbb{R}^{n \times n} \text {,}
\end{equation*}
with
\begin{equation*}
S=\left[\begin{array}{ccccc}
8 & -1 & & & \\
-1 & 8 & -1 & & \\
& \ddots & \ddots & \ddots & \\
& & -1 & 8 & -1 \\
& & & -1 & 8
\end{array}\right]\in \mathbb{R}^{m \times m},n=m^2.
\end{equation*}

\indent The nonlinear terms $(\phi(x),\psi(x))$ of \eqref{VNCP_iterative} are presented on the first column of Tables \ref{tab2}, \ref{tab3} and \ref{tab4}. Accordingly, we let $L_1=L_2=1$ in FPI method \eqref{VNCP_iterative}. Based on this, in Table \ref{tab1}, we list the range of the parameter $\tau$ for the FPI-like methods according to Theorem \ref{theo3.1} and Theorem \ref{theo3.1}. 
In Example 4.2, we set $\mu_1=\mu_2=4$ and $\mu_1=8,\mu_2=4$ to design two different test problems.\\
\indent In Tables \ref{tab2}, \ref{tab3}, and \ref{tab4}, we report the computational results of applying the FPI-like method and the modulus-based matrix splitting method \cite{Fu2024} to Examples 4.1 and 4.2 with different nonlinear terms. The optimal parameter $\alpha$ for FPI-SOR and NMSOR are determined through approximate calculations. We observed that, under the same accuracy level, the FPI-like method required fewer iterations and less computational time than the modulus-based matrix splitting method \cite{Fu2024}, demonstrating superior computational performance. Moreover, the iteration count of the FPI-like method did not significantly increase with the dimensionality $n$, indicating good stability. These results suggest that the FPI-like method can solve more complex problems effectively.
\begin{table}[!ht]
    \caption{The value of $\tau$ for test problems}
    \centering
    \begin{tabular}{lllcc}
    \hline
        test problem & Method & $\tau (\text{Theorem} \ref{theo3.1})$ & $\tau (\text{Theorem} \ref{theo3.2})$ \\ \hline
        Exmaple 4.1 & FPI-J & $(0,1.1667)$ &$(0,1.0667)$\\
        ~ & FPI-GS & $(0,1.1667)$&$(0,1.0833)$\\
        ~ & FPI-SOR & $(0,1.0455)$&$(0,1.0198)$ \\
        \hline 
        Exmaple 4.2($\mu_1=\mu_2=4$) & FPI-J & $(0,1.6031)$ &$(0,1.3367)$\\
        ~ & FPI-GS & $(0,1.6027)$&$(0,1.4021)$\\
        ~ & FPI-SOR & $(0,1.4924)$ &$(0,1.2573)$\\
        \hline 
        Exmaple 4.2($\mu_1=8,\mu_2=4$) & FPI-J & $(0,1.6456)$&$(0,1.3945)$ \\
        ~ & FPI-GS & $(0,1.6452)$&$(0,1.4596)$\\
        ~ & FPI-SOR & $(0,1.5526)$&$(0,1.3123)$ \\
        \hline 
    \end{tabular}
    \label{tab1}
\end{table}
\begin{table}[t]
    \caption{ Numerical comparison of Example 4.1.}
    \centering
    \begin{tabular}{llllllll}
    \hline
        $(\phi(x),\psi(x))$ & Method & \multicolumn{2}{l}{n=1000} & \multicolumn{2}{l}{n=2000} & \multicolumn{2}{l}{n=5000} \\ \cmidrule(lr){3-4} \cmidrule(lr){5-6} \cmidrule(lr){7-8}
        ~ & ~ & IT & Time & IT & Time & IT & Time \\ 
         \hline
        ($|x|,|x|$) & NMJ & 17 & 0.552 & 17 & 2.586 & 17 & 34.923 \\ 
        ~ & NMGS & 15 & 0.498 & 15 & 2.518 & 15 & 31.514 \\ 
        ~ & NMSOR($\alpha=0.8$) & 12 & 0.431 & 12 & 1.852 & 13 & 27.901 \\ 
        ~ & FPI-J($\tau=1.05$) & 10 & 0.176 & 10 & 0.936 & 10 & 13.846 \\ 
        ~ & FPI-GS($\tau=0.95$) & 8 & 0.124 & 8 & 0.609 & 8 & 8.751 \\ 
        ~ & FPI-SOR($\tau=1.01,\alpha=1.05$) & 7 & 0.103 & 7 & 0.555 & 7 & 7.760 \\ \hline
        ($\sin(x),\cos(x)$) & NMJ & 20 & 0.643 & 20 & 3.333 & 20 & 40.992 \\ 
        ~ & NMGS & 14 & 0.443 & 14 & 2.992 & 14 & 29.814 \\ 
        ~ & NMSOR($\alpha=0.9$) & 12 & 0.367 & 12 & 2.430 & 12 & 25.516 \\ 
        ~ & FPI-J($\tau=1.06$) & 15 & 0.328 & 15 & 1.687 & 15 & 24.151 \\ 
        ~ & FPI-GS($\tau=1.04$) & 10 & 0.239 & 10 & 1.285 & 10 & 14.710 \\ 
        ~ & FPI-SOR($\tau=1.01,\alpha=1.05$) & 7 & 0.157 & 7 & 0.973 & 7 & 9.948 \\ \hline
        ($\sin(x),|x|$) & NMJ & 17 & 0.384 & 17 & 2.021 & 17 & 24.443 \\ 
        ~ & NMGS & 12 & 0.314 & 12 & 1.612 & 12 & 19.031 \\ 
        ~ & NMSOR($\alpha=0.9$) & 9 & 0.233 & 9 & 1.114 & 9 & 14.831 \\ 
        ~ & FPI-J ($\tau=1.06$) & 12 & 0.217 & 12 & 1.047 & 12 & 21.040 \\ 
        ~ & FPI-GS($\tau=1.08$) & 7 & 0.131 & 7 & 0.704 & 7 & 16.289 \\ 
        ~ & FPI-SOR($\tau=1.01,\alpha=1.05$) & 6 & 0.114 & 6 & 0.624 & 6 & 8.509 \\ \hline
        ($|x|,\sin(x)$) & NMJ & 16 & 0.525 & 16 & 2.554 & 16 & 33.750 \\ 
        ~ & NMGS & 10 & 0.324 & 10 & 1.689 & 10 & 21.415 \\ 
        ~ & NMSOR($\alpha=0.9$) & 9 & 0.307 & 9 & 1.496 & 10 & 18.267 \\ 
        ~ & FPI-J($\tau=1.05$) & 16 & 0.267 & 16 & 1.440 & 16 & 17.239 \\ 
        ~ & FPI-GS($\tau=1.08$)& 10 & 0.167 & 10 & 0.941 & 10 & 11.304 \\ 
        ~ & FPI-SOR($\tau=1.01,\alpha=1.05$) & 9 & 0.154 & 9 & 0.836 & 9 & 10.150 \\ \hline
        ($|x|,\frac{x}{1+x}$) & NMJ & 22 & 0.513 & 22 & 2.368 & 22 & 30.627 \\ 
        ~ & NMGS & 17 & 0.275 & 17 & 1.631 & 17 & 28.162 \\ 
        ~ & NMSOR($\alpha=0.9$) & 11 & 0.271 & 11 & 1.478 & 11 & 27.801 \\ 
        ~ & FPI-J($\tau=0.97$) & 17 & 0.276 & 17 & 1.191 & 17 & 13.793 \\ 
        ~ & FPI-GS($\tau=1.08$) & 10 & 0.176 & 10 & 0.855 & 10 & 12.504 \\ 
        ~ & FPI-SOR($\tau=1.01,\alpha=1.05$) & 8 & 0.163 & 8 & 0.795 & 8 & 9.399 \\ \hline
        ($\sin(x),\frac{x}{1+x}$) & NMJ & 25 & 0.517 & 25 & 2.664 & 25 & 31.313 \\ 
        ~ & NMGS & 18 & 0.346 & 18 & 1.845 & 18 & 25.414 \\ 
        ~ & NMSOR($\alpha=0.9$) & 14 & 0.302 & 14 & 1.636 & 14 & 24.603 \\ 
        ~ & FPI-J($\tau=0.94$) & 20 & 0.27 & 20 & 1.455 & 20 & 15.502 \\ 
        ~ & FPI-GS($\tau=1.08$) & 14 & 0.177 & 14 & 0.947 & 14 & 12.901 \\ 
        ~ & FPI-SOR($\tau=1.01,\alpha=1.05$) & 11 & 0.174 & 11 & 0.915 & 11 & 10.53 \\ \hline
    \end{tabular}
    \label{tab2}
\end{table}

%% example2 numerical expeirment

\begin{table}[!ht]
    \caption{ Numerical comparison of Example 4.2 for $\mu_1=\mu_2=4$}
    \centering
    \begin{tabular}{llllllll}
    \hline
        $(\phi(x),\psi(x))$ & Method & \multicolumn{2}{l}{n=1600} & \multicolumn{2}{l}{n=3600} & \multicolumn{2}{l}{n=6400} \\ \cmidrule(lr){3-4} \cmidrule(lr){5-6} \cmidrule(lr){7-8}
        ~ & ~ & IT & Time & IT & Time & IT & Time \\ 
         \hline
        ($|x|,|x|$) & NMJ & 17 & 0.862 & 18 & 3.211 & 18 & 33.608 \\ 
        ~ & NMGS & 13 & 0.672 & 13 & 2.221 & 14 & 26.529 \\ 
        ~ & NMSOR($\alpha=0.8$) & 10 & 0.542 & 10 & 1.798 & 10 & 21.021 \\ 
        ~ & FPI-J($\tau=1.10$) & 11 & 0.448 & 11 & 1.512 & 11 & 10.011 \\ 
        ~ & FPI-GS($\tau=1.05$) & 9 & 0.439 & 9 & 1.325 & 9 & 8.393 \\ 
        ~ & FPI-SOR($\tau=1.05,\alpha=1.05$) & 7 & 0.345 & 7 & 0.997 & 7 & 6.546 \\ \hline
        ($\sin(x),\cos(x)$) & NMJ & 27 & 1.796 & 27 & 8.365 & 27 & 57.821 \\ 
        ~ & NMGS & 20 & 1.59 & 20 & 6.298 & 20 & 39.358 \\ 
        ~ & NMSOR($\alpha=0.9$) & 17 & 1.368 & 17 & 5.689 & 17 & 30.706 \\ 
        ~ & FPI-J($\tau=1.10$) & 22 & 1.174 & 22 & 3.224 & 22 & 20.734 \\ 
        ~ & FPI-GS($\tau=0.95$) & 16 & 1.041 & 16 & 2.912 & 16 & 19.531 \\ 
        ~ & FPI-SOR($\tau=1.05,\alpha=1.05$) & 13 & 0.832 & 13 & 2.364 & 13 & 15.697 \\ \hline
        ($\sin(x),|x|$) & NMJ & 17 & 0.921 & 18 & 3.311 & 18 & 24.142 \\ 
        ~ & NMGS & 13 & 0.732 & 13 & 2.192 & 14 & 17.667 \\ 
        ~ & NMSOR($\alpha=0.9$) & 10 & 0.673 & 10 & 2.007 & 11 & 15.847 \\ 
        ~ & FPI-J($\tau=1.10$) & 11 & 0.527 & 11 & 1.765 & 11 & 15.546 \\ 
        ~ & FPI-GS($\tau=1.05$) & 9 & 0.446 & 9 & 1.365 & 9 & 12.874 \\ 
        ~ & FPI-SOR($\tau=1.05,\alpha=1.05$) & 6 & 0.393 & 6 & 1.388 & 6 & 9.235 \\ \hline
        ($|x|,\sin(x)$) & NMJ & 18 & 1.115 & 18 & 6.591 & 18 & 33.272 \\ 
        ~ & NMGS & 14 & 0.859 & 14 & 3.895 & 14 & 28.462 \\ 
        ~ & NMSOR($\alpha=0.9$) & 10 & 0.653 & 11 & 2.35 & 11 & 22.578 \\ 
        ~ & FPI-J($\tau=1.12$) & 12 & 0.631 & 12 & 1.841 & 12 & 11.578 \\ 
        ~ & FPI-GS($\tau=0.95$) & 9 & 0.469 & 10 & 1.308 & 10 & 8.557 \\ 
        ~ & FPI-SOR($\tau=1.10,\alpha=1.05$) & 5 & 0.357 & 5 & 1.273 & 6 & 6.547 \\ \hline
        ($|x|,\frac{x}{1+x}$) & NMJ & 16 & 0.999 & 16 & 2.872 & 16 & 28.293 \\ 
        ~ & NMGS & 12 & 0.744 & 12 & 2.674 & 12 & 25.664 \\ 
        ~ & NMSOR($\alpha=0.9$) & 11 & 0.689 & 11 & 2.15 & 11 & 19.949 \\ 
        ~ & FPI-J($\tau=1.15$) & 9 & 0.368 & 10 & 2.41 & 10 & 15.006 \\ 
        ~ & FPI-GS($\tau=0.95$) & 8 & 0.332 & 8 & 1.921 & 8 & 13.474 \\ 
        ~ & FPI-SOR($\tau=1.05,\alpha=1.05$) & 7 & 0.322 & 7 & 1.452 & 7 & 9.987 \\ \hline
        ($\sin(x),\frac{x}{1+x}$) & NMJ & 17 & 0.914 & 17 & 2.868 & 17 & 16.205 \\ 
        ~ & NMGS & 11 & 0.773 & 11 & 2.705 & 11 & 14.762 \\ 
        ~ & NMSOR($\alpha=0.9$) & 9 & 0.625 & 9 & 2.117 & 9 & 11.556 \\ 
        ~ & FPI-J($\tau=1.12$) & 9 & 0.459 & 10 & 1.78 & 10 & 9.675 \\ 
        ~ & FPI-GS($\tau=0.95$) & 8 & 0.446 & 8 & 1.747 & 8 & 8.542 \\ 
        ~ & FPI-SOR($\tau=1.05,\alpha=1.05$) & 6 & 0.431 & 6 & 1.399 & 6 & 7.372 \\ \hline
    \end{tabular}
    \label{tab3}
\end{table}
\newpage
\begin{table}[!ht]
    \caption{ Numerical comparison of Example 4.2 for $\mu_1=8,\mu_2=4$}
    \centering
    \begin{tabular}{llllllll}
    \hline
        $(\phi(x),\psi(x))$ & Method & \multicolumn{2}{l}{n=1600} & \multicolumn{2}{l}{n=3600} & \multicolumn{2}{l}{n=6400} \\ \cmidrule(lr){3-4} \cmidrule(lr){5-6} \cmidrule(lr){7-8}
        ~ & ~ & IT & Time & IT & Time & IT & Time \\ \hline
         ($|x|,|x|$) & NMJ & 16 & 0.76 & 16 & 2.388 & 16 & 6.817 \\ 
        ~ & NMGS & 12 & 0.558 & 12 & 1.7 & 12 & 4.77 \\ 
        ~ & NMSOR($\alpha=0.8$) & 9 & 0.43 & 9 & 1.432 & 9 & 4.102 \\ 
        ~ & FPI-J($\tau=1.05$) & 9 & 0.313 & 10 & 1.281 & 10 & 4.598 \\ 
        ~ & FPI-GS($\tau=0.95$) & 8 & 0.354 & 8 & 1.073 & 8 & 2.879 \\ 
        ~ & FPI-SOR($\tau=0.95,\alpha=1.05$) & 5 & 0.278 & 5 & 0.937 & 5 & 2.558 \\ \hline
        ($\sin(x),\cos(x)$) & NMJ & 23 & 1.423 & 24 & 4.555 & 24 & 12.33 \\ 
        ~ & NMGS & 18 & 1.105 & 18 & 3.92 & 18 & 10.538 \\ 
        ~ & NMSOR($\alpha=0.9$) & 15 & 1.352 & 15 & 3.29 & 15 & 8.938 \\ 
        ~ & FPI-J($\tau=1.05$) & 17 & 0.915 & 18 & 2.767 & 18 & 6.943 \\ 
        ~ & FPI-GS($\tau=0.95$) & 13 & 0.872 & 13 & 2.553 & 13 & 6.311 \\ 
        ~ & FPI-SOR($\tau=0.95,\alpha=1.05$) & 9 & 0.631 & 10 & 1.843 & 10 & 4.972 \\ \hline
        ($\sin(x),|x|$) & NMJ & 16 & 0.844 & 16 & 2.747 & 16 & 6.913 \\ 
        ~ & NMGS & 12 & 0.599 & 12 & 1.827 & 12 & 5.012 \\ 
        ~ & NMSOR($\alpha=0.9$) & 10 & 0.573 & 10 & 1.826 & 10 & 4.972 \\ 
        ~ & FPI-J($\tau=1.05$) & 9 & 0.422 & 10 & 1.43 & 10 & 3.533 \\ 
        ~ & FPI-GS($\tau=0.95$) & 8 & 0.416 & 8 & 1.222 & 8 & 3.087 \\ 
        ~ & FPI-SOR($\tau=0.95,\alpha=1.05$) & 6 & 0.414 & 6 & 1.16 & 6 & 2.927 \\ \hline
        ($|x|,\sin(x)$) & NMJ & 16 & 0.824 & 16 & 2.507 & 16 & 7.57 \\ 
        ~ & NMGS & 12 & 0.548 & 12 & 1.989 & 13 & 5.908 \\ 
        ~ & NMSOR($\alpha=0.9$) & 11 & 0.638 & 11 & 1.768 & 11 & 5.518 \\ 
        ~ & FPI-J($\tau=1.05$) & 10 & 0.513 & 10 & 1.323 & 10 & 3.617 \\ 
        ~ & FPI-GS($\tau=0.95$) & 8 & 0.416 & 8 & 1.125 & 8 & 3.01 \\ 
        ~ & FPI-SOR($\tau=0.95,\alpha=1.05$) & 6 & 0.39 & 6 & 0.932 & 6 & 2.977 \\ \hline
        ($|x|,\frac{x}{1+x}$) & NMJ & 16 & 0.999 & 16 & 2.872 & 16 & 28.293 \\ 
        ~ & NMGS & 12 & 0.744 & 12 & 2.674 & 12 & 25.664 \\ 
        ~ & NMSOR($\alpha=0.9$) & 11 & 0.689 & 11 & 2.15 & 11 & 19.949 \\ 
        ~ & FPI-J($\tau=1.05$) & 9 & 0.368 & 10 & 2.41 & 10 & 15.006 \\ 
        ~ & FPI-GS($\tau=0.95$) & 8 & 0.332 & 8 & 1.921 & 8 & 13.474 \\ 
        ~ & FPI-SOR($\tau=0.95,\alpha=1.05$) & 7 & 0.322 & 7 & 1.452 & 7 & 9.987 \\ \hline
        ($\sin(x),\frac{x}{1+x}$) & NMJ & 17 & 0.914 & 17 & 2.868 & 17 & 16.205 \\ 
        ~ & NMGS & 11 & 0.773 & 11 & 2.705 & 11 & 14.762 \\ 
        ~ & NMSOR($\alpha=0.9$) & 9 & 0.625 & 9 & 2.117 & 9 & 11.556 \\ 
        ~ & FPI-J($\tau=1.05$) & 9 & 0.459 & 10 & 1.78 & 10 & 9.675 \\ 
        ~ & FPI-GS($\tau=0.95$) & 8 & 0.446 & 8 & 1.747 & 8 & 8.542 \\ 
        ~ & FPI-SOR($\tau=0.95,\alpha=1.05$) & 6 & 0.431 & 6 & 1.399 & 6 & 7.372 \\ \hline
    \end{tabular}
    \label{tab4}
\end{table}
\section{Conclusions}\label{sec5}
In this paper, we propose a novel matrix splitting-based fixed point iteration method (FPI) for solving the VNCP \eqref{VNCP}. By introducing two different iterative error vectors, we derived the convergence conditions for two distinct cases of the FPI method. Notably, we found that the convergence result in Theorem \ref{theo3.2} is a stricter version of that in Theorem \ref{theo3.1}, which is a particularly interesting. Based on the results of Theorem \ref{theo3.2}, we give the estimate of iteration step $k$. Furthermore, by applying two test problems with different nonlinear components, we observe that the FPI method outperforms the modulus-based matrix splitting method \cite{Fu2024} in terms of both iterations(iter) and computation time(CPU). In future work, we will focus on theoretically identifying the optimal parameter $\tau^{\text{opt}}$.
\newpage
\section*{Data availability} {\footnotesize
Data sharing is not applicable to this article as no datasets were generated or analyzed during the current study.
\section*{Conflict of interest}
We declare that we have no Conflict of interest in this paper.
\bmhead{Acknowledgements}
This work was supported by the Natural Science Basic Research Program of Shaanxi (Program No. 2024JC-YBMS-026).

\bibliography{sn-bibliography}

\end{document}